\documentclass{svproc}

\usepackage{graphicx}
\usepackage{amsmath,amssymb}
\usepackage{booktabs}
\usepackage{multirow}
\usepackage{url}

\usepackage[hidelinks]{hyperref}

\graphicspath{{figures/}}

\begin{document}
\mainmatter

\title{Constrained Spatial Pricing of On-Street Parking with Bayesian Demand Calibration}
\titlerunning{Constrained Spatial Parking Pricing}
\author{
Ananya Kale\inst{1} \and
Mohit Apte\inst{2}
}

\authorrunning{Ananya Kale and Mohit Apte}
\tocauthor{Ananya Kale, Mohit Apte}

\institute{
MKSSS's Cummins College of Engineering for Women, India\\
\email{ananya.kale@cumminscollege.in}
\and
University of Chicago, Chicago, IL, USA\\
\email{mohitapte@uchicago.edu}
}
\maketitle

\begin{abstract}
Occupancy-targeted curb pricing, including San Francisco's SFpark pilot, often relies on
local threshold rules that do not account for spatial substitution or uncertainty in demand
response. We estimate parking-price elasticities from disaggregated SFpark data using a
hierarchical Bayesian demand model and use the resulting posterior distribution in a
spatially coupled constrained pricing model. Partial pooling across districts gives a
posterior mean elasticity of $-0.241$ with a $95\%$ credible interval of
$[-0.325,-0.176]$. The pricing problem is formulated as a constrained nonlinear program
with spatial coupling, temporal price-change limits, neighbor price gaps, and a soft
revenue floor, and is solved with Ipopt. Expected-loss and CVaR objectives are evaluated
over joint posterior draws rather than a fixed elasticity. On held-out posterior scenarios,
the posterior-based CVaR policy has a lower objective than the SFpark threshold rule and
historical tariff with posterior probability $1.00$, and than a literature-calibrated robust
policy with probability $0.85$. Epsilon-constraint frontiers show that posterior
calibration changes the efficient pricing set. Results from 60 district-windows across
Fillmore, Mission, and Marina also report the associated tradeoffs in occupancy-band
performance, revenue, and spatial disparity.

\end{abstract}
\section{Introduction}
Chronic over- and under-occupancy of curb space is a durable urban operations problem.
When prices are fixed, popular blocks saturate while nearby capacity sits idle; drivers
cruise for vacancies, adding congestion~\cite{arnott2006,shoup2006}. SFpark adjusted
on-street rates toward an occupancy target band using sensor data and periodic rate
changes~\cite{pierce2013,millard2014}. Operational practice, however, remains largely
local: each block's rate responds to its own occupancy history. Local threshold rules do
not coordinate neighboring prices, even though demand substitutes across space, and they
do not hedge uncertainty in how occupancy responds to price.

Two issues complicate this formulation. First, SFpark raised prices where occupancy was already
high, so panel associations between price and occupancy are endogenous~\cite{pierce2013};
in our own attempts, both simple regressions and flexible double-ML residualizations on
this panel produced unstable or wrong-signed elasticities. Second, meaningful spatial
coupling requires a neighbor relation $W$; the public SFpark files used here contain street
and block-number fields but not latitude/longitude.

We address the first difficulty with a Bayesian calibration rather than a design-based
identification strategy. The pricing decision process is reconstructed at the grain at
which rates were actually set (block $\times$ revision wave $\times$ time-of-day band
$\times$ day type), and a hierarchical Bayesian demand model combines the observed SFpark
likelihood with prior information from the parking-elasticity literature and partial
pooling across districts, under an economically motivated negative-support
parameterization. We do not claim nonparametric causal identification. We use the resulting
posterior distribution of short-run price response directly in the pricing model, entering
the optimization as full posterior draws rather than a point summary.

\paragraph{Research question.}
Given a posterior distribution over the demand plant's price response, how should a
neighborhood set block-hour prices that trade occupancy-band performance against revenue
and spatial price consistency, and with what posterior probability do the resulting
policies outperform threshold rules and literature-calibrated robust benchmarks?

\paragraph{Contributions.}
\begin{enumerate}
\item A hierarchical Bayesian behavioral calibration of parking-price response estimated on
      disaggregated SFpark pricing-band observations, with a sign-constrained
      parameterization $\varepsilon_d=-\exp(\eta_d)$, district partial pooling, explicit
      prior-to-posterior updating, and posterior predictive validation.
\item A fully specified constrained NLP with independent and spatially coupled variants,
      solved by Ipopt with multistart, and mapped to meter increments.
\item Posterior-aware optimization: expected-loss and CVaR objectives over joint posterior
      parameter draws, preserving district heterogeneity and parameter correlations, plus
      literature-mapped robust counterparts as benchmarks.
\item Probabilistic policy statements of the form
      $P(J_{\mathrm{Bayes}}<J_{\mathrm{threshold}}\mid\mathcal D)$ that connect estimation
      uncertainty directly to policy choice.
\item Epsilon-constraint Pareto frontiers between band performance, revenue, and spatial
      disparity under three plant calibrations, showing that the estimated posterior
      changes the efficient decision set, together with multi-construction spatial graphs,
      a multi-district calendar evaluation (Fillmore, Mission, Marina; 60
      district-windows), and a small-instance optimality benchmark for Ipopt.
\end{enumerate}

\section{Related work}
\paragraph{Cruising, curb scarcity, and SFpark.}
Economic models of downtown parking link underpriced curb space to cruising and
congestion~\cite{arnott2006,arnott2009,shoup2006}. SFpark is the leading empirical
implementation of occupancy-targeted pricing; Pierce and Shoup~\cite{pierce2013} document
time-of-day and location-varying rate changes and report a mean on-street occupancy
elasticity near $-0.4$ from first-year adjustments, with large heterogeneity. Millard-Ball
et al.~\cite{millard2014} assess performance-parking impacts and caution that short-run
elasticities may be weaker than headline averages. Chatman and Manville~\cite{chatman2014}
and related planning evaluations emphasize institutional and placard-abuse constraints that
attenuate price response.

\paragraph{Parking elasticities.}
Broader reviews place typical parking price elasticities roughly in
$[-0.1,-0.6]$, with a central value near $-0.3$~\cite{vaca2005,litman2012}. These ranges
enter our study in two ways: they define the ambiguity set for a literature-calibrated
robust benchmark, and they inform the prior of a hierarchical Bayesian model that is then
updated with the disaggregated SFpark panel. Constrained parameterizations of the kind we
use to encode a sign restriction are common in applied Bayesian
modeling~\cite{gelman2013,salvatier2016}.

\paragraph{Dynamic pricing and spatial substitution.}
Subsequent work studies occupancy-driven or predictive pricing, including stochastic control
and dynamic programming formulations~\cite{qian2014,qian2015} and performance-based
management~\cite{mackowski2015}; these often assume a known response surface or optimize
predictive accuracy without the operational constraints (temporal caps, neighbor gaps,
revenue floors, meter increments) that agencies must satisfy. Network and equilibrium
models of parking search~\cite{leurent2012,boyles2015,zakharenko2016} provide richer
spatial structure but are difficult to operationalize as day-to-day meter-setting tools;
our coupling is deliberately reduced-form, with neighbor-weighted prices and lags entering
an embeddable occupancy map and neighbor price gaps entering as hard constraints.

\paragraph{Robust and risk-sensitive pricing.}
When demand parameters are uncertain, robust and distributionally aware formulations are
standard~\cite{ben2009,bertsimas2011}. Conditional value-at-risk (CVaR) provides a
tractable tail objective~\cite{rockafellar2000,rockafellar2002}. Applications to
transportation pricing under demand ambiguity motivate scenario and epigraph
forms~\cite{yin2009}.

\paragraph{Gap.}
Existing work either assumes known price response, models spatial interactions through
richer equilibria that are hard to operationalize as meter policies, or uses local pricing
rules. This paper studies \emph{implementable} spatially coupled pricing when elasticity
itself is ambiguous, with explicit NLP/robust formulations, multi-graph sensitivity, and
multi-district calendar evidence.

\section{Demand plant and Bayesian elasticity calibration}
\label{sec:ambiguity}
\subsection{Neighborhood instance}
Let $\mathcal{B}$ be curb blocks in one district and
$\mathcal{T}=\{0,\ldots,T-1\}$ a horizon of consecutive operating hours. Decision variables
are prices $p_{b,t}$. Occupancy $o_{b,t}\in(0,1)$ is determined by the plant below.
Capacities $c_b$ and historical prices $p^{\mathrm{hist}}_{b,t}$ are data. A row-normalized
neighbor matrix $W$ encodes spatial adjacency (Section~\ref{sec:graphs}).

\subsection{Fractional-logit occupancy plant}
\begin{equation}
\label{eq:occ}
o_{b,t}=\sigma(\eta_{b,t}),\qquad
\sigma(z)=\bigl(1+e^{-z}\bigr)^{-1},
\end{equation}
equivalently the implicit form $o_{b,t}\bigl(1+\exp(-\eta_{b,t})\bigr)=1$ used in the NLP.
The linear predictor is
\begin{align}
\label{eq:eta}
\eta_{b,t}
&=
\bar\eta_{b,t}
+\beta_p\,z(p_{b,t})
+\beta_n\,z\bigl((Wp_{\cdot,t})_b\bigr)
\nonumber\\
&\quad
+\beta_\ell\,\operatorname{logit}(o_{b,t-1})
+\beta_{n\ell}\,(Wo_{\cdot,t-1})_b,
\end{align}
where $\bar\eta_{b,t}$ collects block, hour, and day-of-week effects, $z(\cdot)$ denotes
standardization with sample mean $\mu_p$ and standard deviation $s_p$, and $o_{b,-1}$ is an
initial occupancy. The \emph{independent} model sets $\beta_n=\beta_{n\ell}=0$ and drops
neighbor price-gap constraints.

\subsection{Mapping published elasticities to $\beta_p$}
Local occupancy elasticity satisfies
$\varepsilon=(1-o)\,\beta_p\,p/s_p$, hence
\begin{equation}
\label{eq:map}
\beta_p=\frac{\varepsilon\,s_p}{(1-o)\,p}.
\end{equation}
Evaluating at $(\mu_p,s_p,o)=(2.73,0.937,0.70)$ from the fitted plant maps the review range
$\varepsilon\in[-0.6,-0.1]$ to $\beta_p\in[-0.686,-0.114]$, with the central
$\varepsilon=-0.30$ mapping to $\beta_p=-0.343$.

The literature-calibrated benchmark plant uses $\beta_p=-0.343$ ($\varepsilon=-0.30$), and
its robust scenarios draw $\beta_p$ uniformly from $[-0.686,-0.114]$, with
neighbor-coefficient multipliers and initial-occupancy shifts as secondary perturbations.
The posterior-based plant instead maps each joint posterior draw of the district-specific
elasticity into $\beta_p$ via~\eqref{eq:map}, as described next.

\subsection{Hierarchical Bayesian calibration on pricing-band cells}
\label{sec:bayes}
SFpark set rates separately by time-of-day band and weekday/weekend, so the decision panel
is reconstructed at that grain: one cell per block $\times$ revision wave $\times$ day type
$\times$ time-of-day band, with $D_i=\Delta\log p_i$ (posted rate change) and
$Y_i=\Delta\log o_i$ (14-day post- vs.\ pre-revision occupancy). Aggregating to block-wave
cells instead lets opposing within-block band changes cancel; when we fit the same model on
that coarser panel, both the price signal and posterior predictive adequacy were lost.
On the band-level panel
($n{=}10{,}659$ cells; $6{,}721$ with nonzero price changes) we estimate
\begin{align}
\label{eq:b6}
Y_i &= \alpha + \varepsilon_{d(i)}\,D_i + \gamma\,\tilde o_i + u_{w(i)} + v_{g(i)}
      + e_i,\qquad e_i\sim\mathcal N(0,\sigma^2),\\
\varepsilon_d &= -\exp(\eta_d),\qquad
\eta_d\sim\mathcal N(\mu_\eta,\sigma_\eta^2),
\nonumber
\end{align}
where $\tilde o_i$ is a leakage-safe mean-reversion control (pre-revision occupancy
deviation from the cell's own earlier-wave history), $u_w$ and $v_g$ are wave and band
random effects, and $d(i)$ indexes districts. The transformation
$\varepsilon_d=-\exp(\eta_d)$ encodes the economic restriction that, holding the modeled
demand state fixed, raising price does not raise demand; it is a structural modeling choice
analogous to positive-support variance parameters. The hyperprior on $\mu_\eta$ is set so
the pooled elasticity has an informative literature prior,
$\varepsilon\sim\mathcal N(-0.30,0.05)$ truncated to $\varepsilon<0$. Posteriors are sampled
with NUTS in PyMC~\cite{salvatier2016,hoffman2014} ($\max\hat R=1.04$, zero divergences;
rank plots for the pooled location and hierarchical scale in Fig.~\ref{fig:rank}, full
trace and energy diagnostics in supplementary material).

\begin{figure}[!t]
\centering
\includegraphics[width=0.92\textwidth]{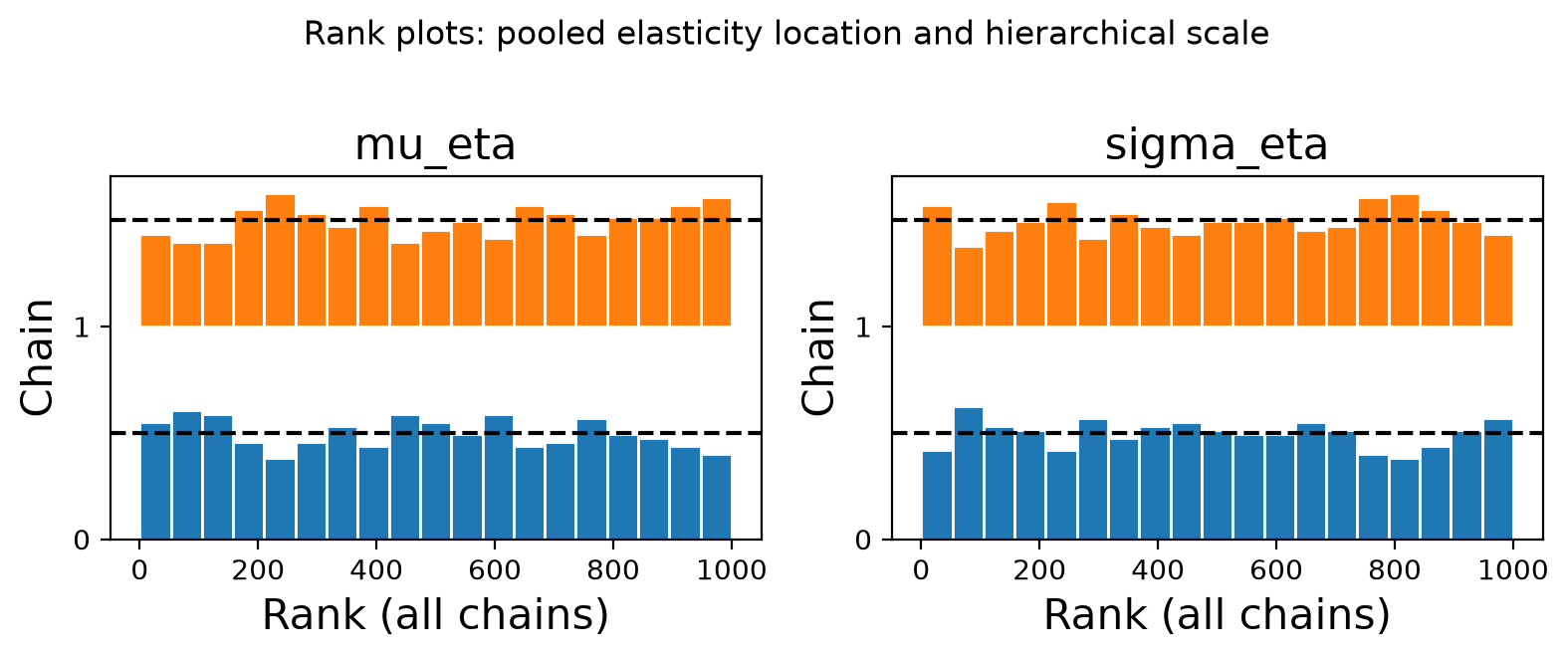}
\caption{Rank plots for the pooled elasticity location $\mu_\eta$ and hierarchical scale
$\sigma_\eta$; near-uniform ranks across chains indicate good mixing.}
\label{fig:rank}
\end{figure}

\begin{figure}[!t]
\centering
\includegraphics[width=0.78\textwidth]{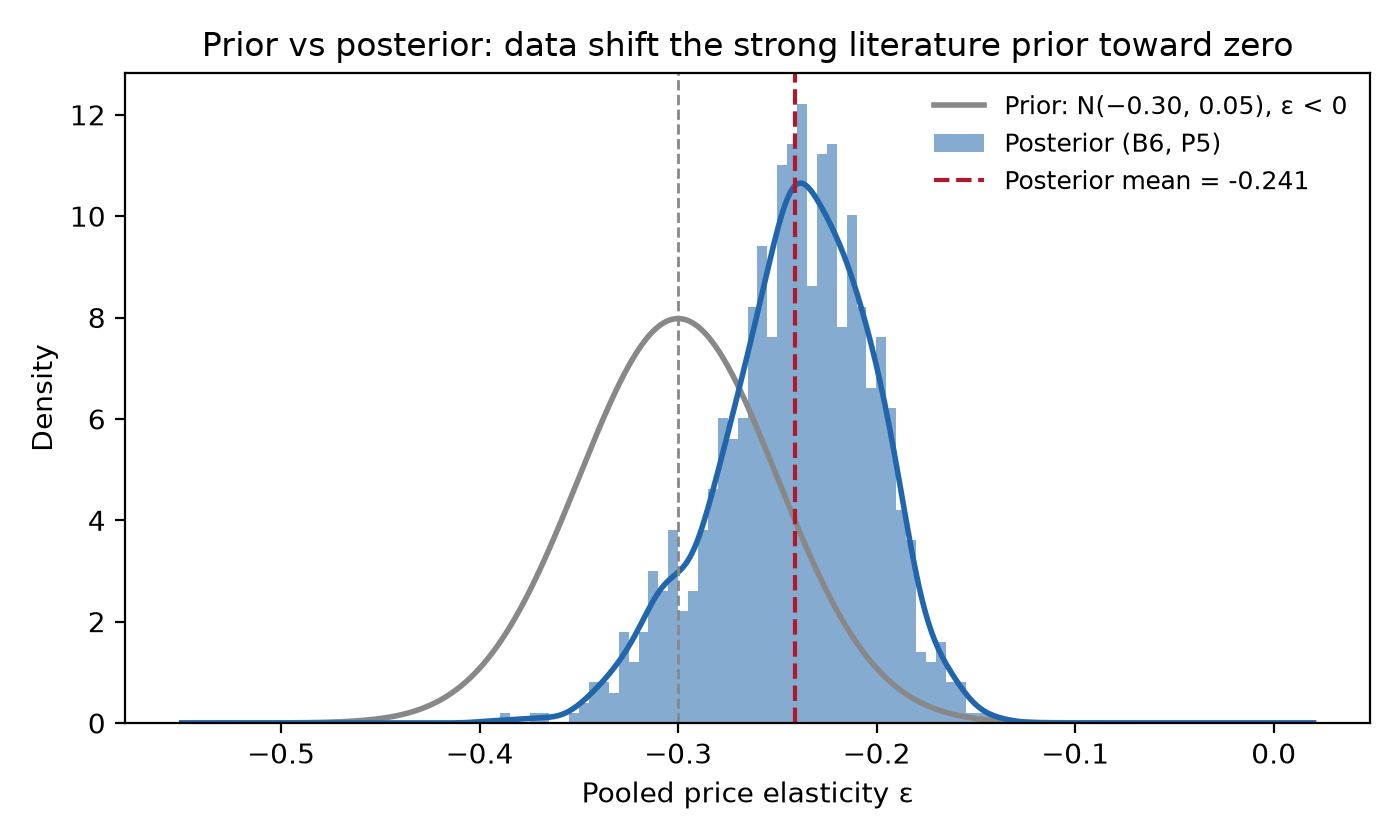}
\caption{Prior versus posterior for the pooled elasticity. The strong literature prior
$\mathcal N(-0.30,0.05)$ (truncated to $\varepsilon<0$) is updated by the band-level SFpark
panel to a posterior with mean $-0.241$.}
\label{fig:priorpost}
\end{figure}

\begin{figure}[!t]
\centering
\includegraphics[width=0.78\textwidth]{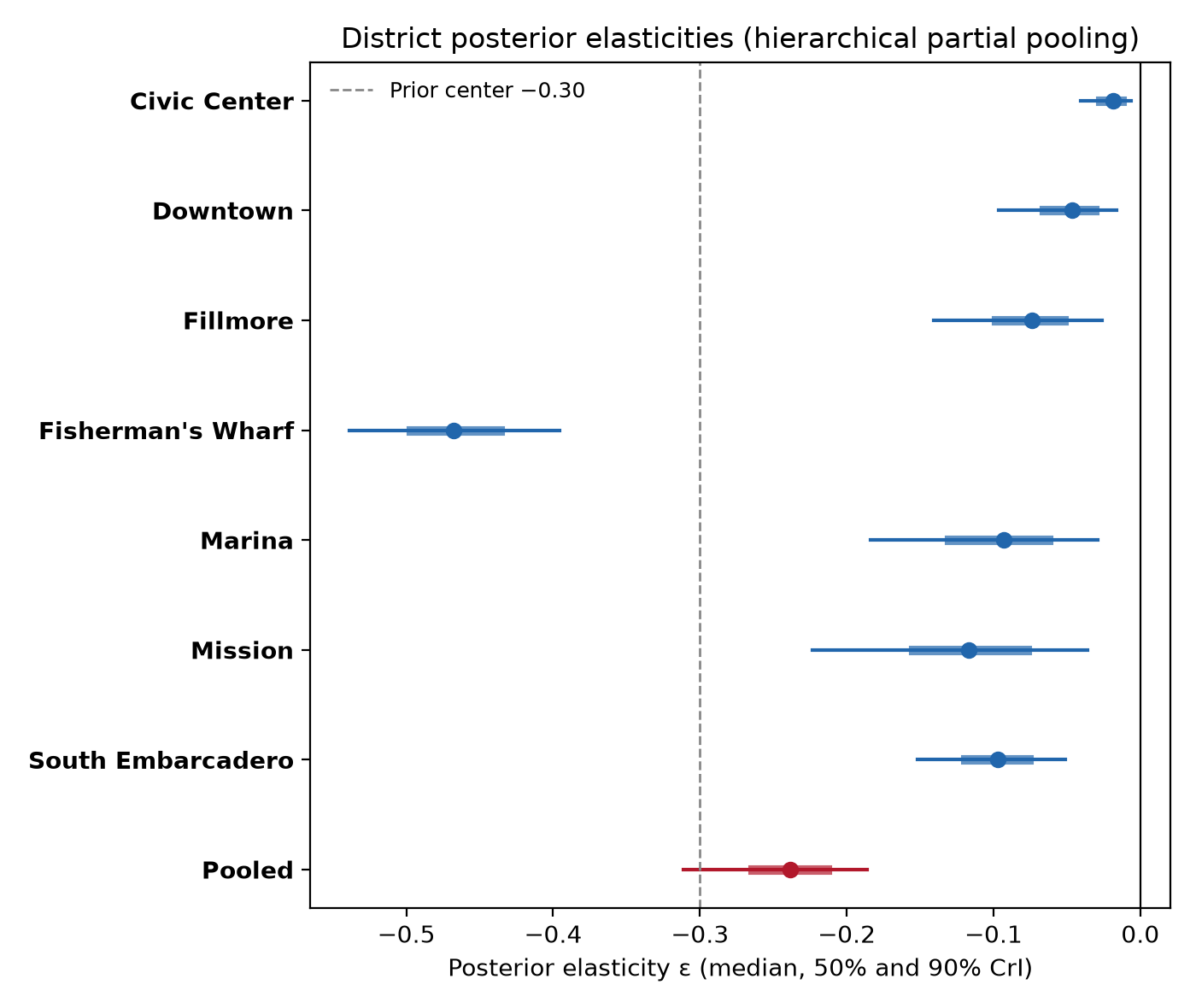}
\caption{District-specific posterior elasticities (medians with 50\% and 90\% credible
intervals) under hierarchical partial pooling, with the pooled posterior at bottom.}
\label{fig:forest}
\end{figure}

The pooled posterior is
\begin{equation}
\label{eq:posterior}
\varepsilon=-0.241,\qquad 95\%\ \mathrm{CrI}=[-0.325,-0.176],
\end{equation}
which shifts the prior center of $-0.30$ noticeably toward zero
(Fig.~\ref{fig:priorpost}); the posterior does not simply reproduce the prior, so the data
are informative. A companion fit with a sign-only, weakly informative prior also produces a
negative posterior ($-0.069$, $[-0.131,-0.028]$), which indicates that the likelihood on
its own favors a negative response. District posteriors differ substantially
(Fig.~\ref{fig:forest}): tourist-heavy Fisherman's Wharf is the most elastic district
(median near $-0.47$), while Civic Center and Downtown are close to inelastic. Empirical
coverage of $80\%$ posterior predictive intervals is $0.88$, and leave-one-wave-out refits
move the pooled mean by at most a few hundredths. Throughout the paper we refer to
\eqref{eq:posterior} as a posterior elasticity estimate from a prior-informed hierarchical
calibration; we do not describe it as causally identified, and no statistical model is free
of assumptions.

\section{Constrained NLP formulation}
\label{sec:nlp}
\subsection{Auxiliary variables}
Target band $[\underline o,\overline o]=[0.6,0.8]$. Define
\begin{align}
\mathrm{over}_{b,t} &\ge o_{b,t}-\overline o,
&
\mathrm{over}_{b,t}&\ge 0,
\label{eq:over}\\
\mathrm{under}_{b,t} &\ge \underline o-o_{b,t},
&
\mathrm{under}_{b,t}&\ge 0,
\label{eq:under}
\end{align}
and a cruising / search-pressure proxy
\begin{equation}
\label{eq:cruise}
s(o)=\frac{1}{1.02-o}-\frac{1}{1.02},
\end{equation}
which rises sharply as occupancy approaches one. Let $r\ge 0$ be a revenue shortfall
auxiliary. Temporal increments are
\begin{equation}
\label{eq:dprice}
\Delta p_{b,0}=p_{b,0}-p^{\mathrm{hist}}_{b,0},
\qquad
\Delta p_{b,t}=p_{b,t}-p_{b,t-1}\ (t\ge 1).
\end{equation}

\subsection{Objective}
With historical-scenario normalization constants $N_{\cdot}>0$ (so that no term dominates by
units) and nonnegative weights $w_{\cdot}$, the deterministic program minimizes
\begin{align}
\label{eq:obj}
J(p,o)
&=
w_o\frac{1}{|\mathcal{B}|T}\sum_{b,t}\mathrm{over}_{b,t}^2\Big/N_o
+
w_u\frac{1}{|\mathcal{B}|T}\sum_{b,t}\mathrm{under}_{b,t}^2\Big/N_u
\nonumber\\
&\quad
+
w_s\frac{1}{|\mathcal{B}|T}\sum_{b,t}s(o_{b,t})\Big/N_s
+
w_\tau\frac{1}{|\mathcal{B}|T}\sum_{b,t}(\Delta p_{b,t})^2\Big/N_\tau
\nonumber\\
&\quad
+
w_x\frac{1}{|\mathcal{B}|T}\sum_{t}\sum_{i,j}W_{ij}(p_{i,t}-p_{j,t})^2\Big/N_x
+
w_r\,r\Big/N_r.
\end{align}
Default weights are
$(w_o,w_u,w_s,w_\tau,w_x,w_r)=(1.0,0.5,1.0,0.25,0.25,0.5)$.
Because~\eqref{eq:obj} is a weighted sum, band share is not maximized as a single KPI;
Section~\ref{sec:frontier} treats that mismatch as the object of study.

\subsection{Constraints}
\begin{align}
p_{\min} &\le p_{b,t} \le p_{\max},
\label{eq:bounds}\\
|\Delta p_{b,t}| &\le \Delta^{\max},
\label{eq:temp}\\
|p_{i,t}-p_{j,t}| &\le \Gamma
\quad\text{whenever }W_{ij}>0\text{ (coupled model)},
\label{eq:gap}\\
\sum_{b,t} p_{b,t}\,o_{b,t}\,c_b + r
&\ge
\rho\sum_{b,t} p^{\mathrm{hist}}_{b,t}\,o^{\mathrm{hist}}_{b,t}\,c_b,
\label{eq:rev}
\end{align}
together with the occupancy equations~\eqref{eq:occ}--\eqref{eq:eta} and auxiliaries
\eqref{eq:over}--\eqref{eq:under}. In the study,
$[p_{\min},p_{\max}]=[0.25,8]$, $\Delta^{\max}=\$1$, $\Gamma=\$2$, and $\rho=0.9$.
The independent model omits~\eqref{eq:gap} and neighbor terms in~\eqref{eq:eta}.

\subsection{Solution method and implementability}
The NLP is implemented in Pyomo and solved with Ipopt~\cite{wachter2006}. Multistart uses
historical, threshold, uniform, random, and perturbed-best initial prices; we report the
best local solution and do not claim global optimality. After each solve, occupancy and
objective components are recomputed by an independent NumPy simulator.

Continuous prices are rounded to a \$0.25 grid, repaired to restore~\eqref{eq:bounds},
\eqref{eq:temp}, and~\eqref{eq:gap}, then improved by coordinate-wise discrete local search
over $\{\pm 0.25,0\}$ moves that preserve feasibility.

\section{Robust counterparts}
\label{sec:robust}
\subsection{Scenarios}
Uncertainty is represented by $K$ scenarios. Scenario $k$ draws an own-price coefficient
$\beta_p^{(k)}$ from the literature interval in Section~\ref{sec:ambiguity}, multiplies the
neighbor-price coefficient by $\nu^{(k)}\sim\mathrm{Unif}(1\pm 0.5)$, and shifts initial
occupancy by $\delta^{(k)}\sim\mathrm{Unif}([-0.08,0.08])$ (stress scenarios use larger
shifts). Let $J_k(p)$ denote objective~\eqref{eq:obj} under the occupancy path induced by
scenario $k$.

\paragraph{Shared vs.\ scenario-specific quantities.}
Prices $p=\{p_{b,t}\}$ are \emph{shared} across scenarios (one implementable tariff).
Occupancy paths $o^{(k)}$, auxiliaries $(\mathrm{over}^{(k)},\mathrm{under}^{(k)},r^{(k)})$,
and therefore $J_k(p)$ are \emph{scenario-specific}. Constraints~\eqref{eq:bounds}--\eqref{eq:gap}
apply to the shared $p$; the revenue floor~\eqref{eq:rev} is imposed per scenario with its
own shortfall $r^{(k)}$.

\subsection{Aggregation modes}
\begin{align}
\text{Expected:}\quad
&\min_p\ \frac{1}{K}\sum_{k=1}^{K} J_k(p).
\label{eq:exp}\\[0.4em]
\text{Worst-case:}\quad
&\min_{p,z}\ z
\quad\text{s.t.}\quad
z\ge J_k(p)\ \forall k.
\label{eq:wc}\\[0.4em]
\text{CVaR$_\alpha$:}\quad
&\min_{p,\zeta,u}\
\zeta+\frac{1}{(1-\alpha)K}\sum_{k=1}^{K} u_k
\label{eq:cvar}\\
&\quad\text{s.t.}\quad
u_k\ge J_k(p)-\zeta,\quad u_k\ge 0,
\quad k=1,\ldots,K.
\nonumber
\end{align}
We use $\alpha=0.8$ and $K=8$ optimization scenarios, with held-out validation ($K=8$) and
stress sets for evaluation. The literature benchmark's robust mode is CVaR.

\subsection{Posterior-aware counterparts}
\label{sec:postopt}
For the posterior-based policies, scenarios are joint draws
$\theta^{(k)}\sim p(\theta\mid\mathcal D)$ from the fitted hierarchical posterior of
Section~\ref{sec:bayes} rather than uniform draws from a literature interval. Each draw's
district-specific elasticity is mapped into $\beta_p^{(k)}$ via~\eqref{eq:map} using the
posterior component of the district being priced, so both district heterogeneity and
parameter correlations are retained in the scenario set. We solve
\begin{equation}
\label{eq:postobj}
\min_p\ \mathbb E_{\theta\mid\mathcal D}\bigl[J(p,\theta)\bigr]
\qquad\text{and}\qquad
\min_p\ \mathrm{CVaR}_\alpha\bigl(J(p,\theta)\bigr),
\end{equation}
using the same shared-price scenario program as~\eqref{eq:exp}--\eqref{eq:cvar}. Because the
scenario distribution is now a posterior, evaluations of any fixed policy across held-out
posterior draws yield direct probabilistic statements such as
$P\bigl(J(p_{\mathrm{Bayes}},\theta)<J(p_{\mathrm{threshold}},\theta)\mid\mathcal D\bigr)$.

\section{Spatial graphs}
\label{sec:graphs}
Absent coordinates in the ingested SFpark files, we do not claim a true geographic network.
We compare four constructions using observed street/block fields and a schematic layout:
street adjacency (same $\mathtt{street\_name}$, block-number gap at most 4, distance-decay
weights, row-normalized), street binary (same support, uniform weights), schematic $k$-NN
($k=6$ in a street/block schematic embedding), and schematic distance (edges within a
threshold in the same embedding). Coupled-price conclusions are required to persist across
constructions (Section~\ref{sec:results}).

\section{Experimental design}
\label{sec:design}
\paragraph{Data and instances.}
SFpark hourly occupancy and rates (2011--2013)~\cite{sfmta2014}. Primary NLP instance:
Fillmore, 12 highest-volume high-quality blocks, $T=6$ operating hours. Multi-district
calendar evaluation: Fillmore, Mission, and Marina, each with 12 blocks and 20 evenly spaced
weekday windows with full block coverage (60 district-windows). Fillmore also has a
24-window comparison including robust CVaR.

\paragraph{Policies.}
Historical prices; a threshold rule that steps rates toward the band; independent NLP;
spatially coupled NLP; literature-robust coupled CVaR (Fillmore 24-window study); and the
posterior-based expected-loss and CVaR policies of Section~\ref{sec:postopt}, evaluated on
held-out posterior draws.

\paragraph{Frontiers and optimality.}
Pareto frontiers between band performance, revenue, and disparity are generated with the
epsilon-constraint method of Section~\ref{sec:frontier} under three plant calibrations. On
$3\times 3$, $4\times 3$, and $4\times 4$ instances we compare Ipopt's best multistart
objective to the best feasible point among hundreds of increment-grid candidates.

\section{Results}
\label{sec:results}
\subsection{Elasticity sensitivity}
Across the mapped literature grid, the coupled model retains a revenue advantage and a
neighbor-disparity advantage over the threshold rule for every tested
$\varepsilon\in[-0.6,-0.05]$; there is no revenue break-even inside the review range, and
band share remains threshold-led under all tested elasticities, consistent with
objective~\eqref{eq:obj}.

\subsection{Single-scenario Fillmore profile}
With the central plant, Ipopt succeeds on $7/7$ multistarts for independent and coupled
models, with multistart objectives agreeing to relative spread $<10^{-8}$ and external
occupancy verification to $4.5\times 10^{-8}$. Table~\ref{tab:indnet} shows the KPI split
that motivates the frontier analysis.

\begin{table}[!t]
\caption{Fillmore representative scenario ($12\times 6$), central plant $\beta_p=-0.343$.}
\label{tab:indnet}
\centering
\scriptsize
\begin{tabular}{@{}lrrrr@{}}
\toprule
Policy & Band share & Revenue & Disparity & Cruising \\
\midrule
Historical  & 0.792 & 2231 & 0.079 & 1.780 \\
Threshold   & 0.944 & 1945 & 0.088 & 2.053 \\
Independent & 0.875 & 2546 & 0.014 & 1.408 \\
Coupled     & 0.083 & 2560 & 0.012 & 1.391 \\
\bottomrule
\end{tabular}
\end{table}

\subsection{Pareto frontiers by the epsilon-constraint method}
\label{sec:frontier}
A weighted-sum sweep of~\eqref{eq:obj} can recover efficient points, but it may generate
dominated solutions, miss nonconvex parts of the frontier, and depends on normalization.
We therefore construct the band--revenue--disparity frontier directly with the
epsilon-constraint method. For a grid of revenue floors
$R_{\min}\in\{0.90,0.95,\ldots,1.30\}\times R_{\mathrm{hist}}$ and three neighbor-disparity
caps $D_{\max}$ (strict, medium, loose; $0.5\times$, $2\times$, and $8\times$ the
historical disparity), we solve
\begin{equation}
\label{eq:epscon}
\min_p\ \mathbb E_k\bigl[\mathrm{BandViolation}(p,\theta_k)\bigr]
\quad\text{s.t.}\quad
\mathbb E_k\bigl[R(p,\theta_k)\bigr]\ge R_{\min},
\qquad
D(p)\le D_{\max},
\end{equation}
retaining all operational constraints
\eqref{eq:bounds}--\eqref{eq:gap}. Band violation is the mean per block-hour of
$\mathrm{over}+\mathrm{under}$. The sweep is run under three calibrations of the scenario
set $\{\theta_k\}$: the deterministic central plant ($\varepsilon=-0.30$), the
literature-robust scenario set, and joint draws from the B6 posterior. Solutions are
re-evaluated on held-out scenario sets, and dominated points are removed by a
three-objective nondominance filter (81 solves; each nondominated point is an optimized
tariff).

\begin{figure}[!t]
\centering
\includegraphics[width=\textwidth]{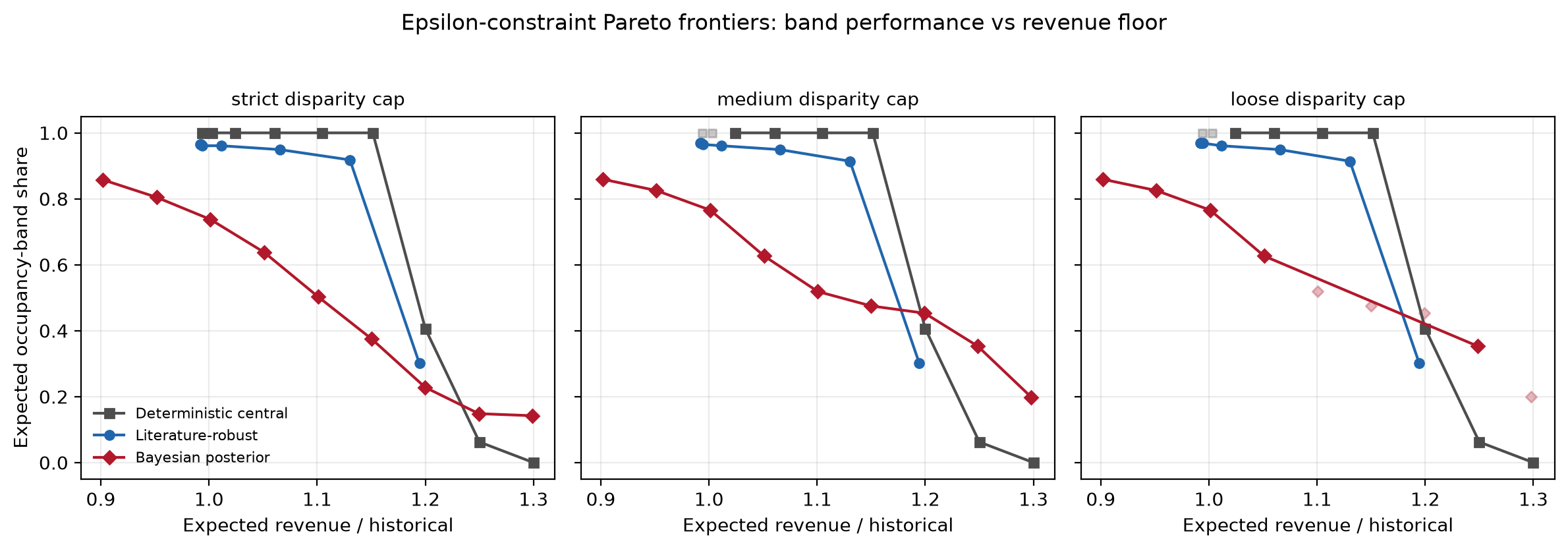}
\caption{Epsilon-constraint Pareto frontiers between expected occupancy-band share and the
expected revenue floor (relative to historical), under strict, medium, and loose
neighbor-disparity caps. Faded markers are dominated points; missing high-revenue points
were certified locally infeasible.}
\label{fig:pareto}
\end{figure}

Figure~\ref{fig:pareto} shows the result. Under the deterministic and literature
calibrations the frontier is nearly flat: full band compliance is compatible with revenue
floors up to about $1.15\times$ historical, after which the frontier collapses, with knees
at $1.20\times$ (deterministic) and $1.15\times$ (literature); the literature frontier is
infeasible beyond $1.20\times$ because its most elastic scenarios cannot meet the floor.
The posterior frontier is qualitatively different. Because the Fillmore posterior response
is weaker than the assumed $-0.30$, band share declines steadily from the first revenue
floor, the knee arrives earlier (near $1.05$--$1.10\times$ historical), and the local
marginal tradeoff along adjacent nondominated points is roughly $3.7$--$4.7$ revenue-proxy
units per percentage point of band share near historical revenue, steepening beyond the
knee. Beyond about $1.2\times$ historical the posterior frontier crosses above the other
two: the assumed-elasticity calibrations overstate what high revenue floors cost in band
performance at that range, while understating the cost at moderate floors. The posterior calibration therefore changes the efficient set, not only the selected tariff. The
disparity cap is nearly free for the posterior and literature
frontiers (the strict and loose curves almost coincide), which is consistent with the
neighbor-gap constraints already limiting price dispersion.

\subsection{Graph construction sensitivity}
Under all four $W$ constructions of Section~\ref{sec:graphs}, coupled prices have weakly
higher revenue and weakly lower disparity than independent prices on the same plant.
Street-adjacency graphs, which use only observed street names and block numbers, produce
the same qualitative ordering as schematic $k$-NN, so the coupled-price conclusions do not
depend on the schematic embedding.

\subsection{Robust modes}
On held-out validation scenarios drawn from the literature $\beta_p$ interval, CVaR and
worst-case improve mean and worst objectives relative to the deterministic coupled solve
(validation worst $0.900\to 0.869$; mean $0.834\to 0.829$). Stress worst objectives likewise
improve ($1.195\to 1.170$). Robust aggregation improves tails rather than mean band share.

\subsection{Posterior-aware policies and probabilistic policy claims}
\label{sec:postresults}
\begin{figure}[!t]
\centering
\includegraphics[width=\textwidth]{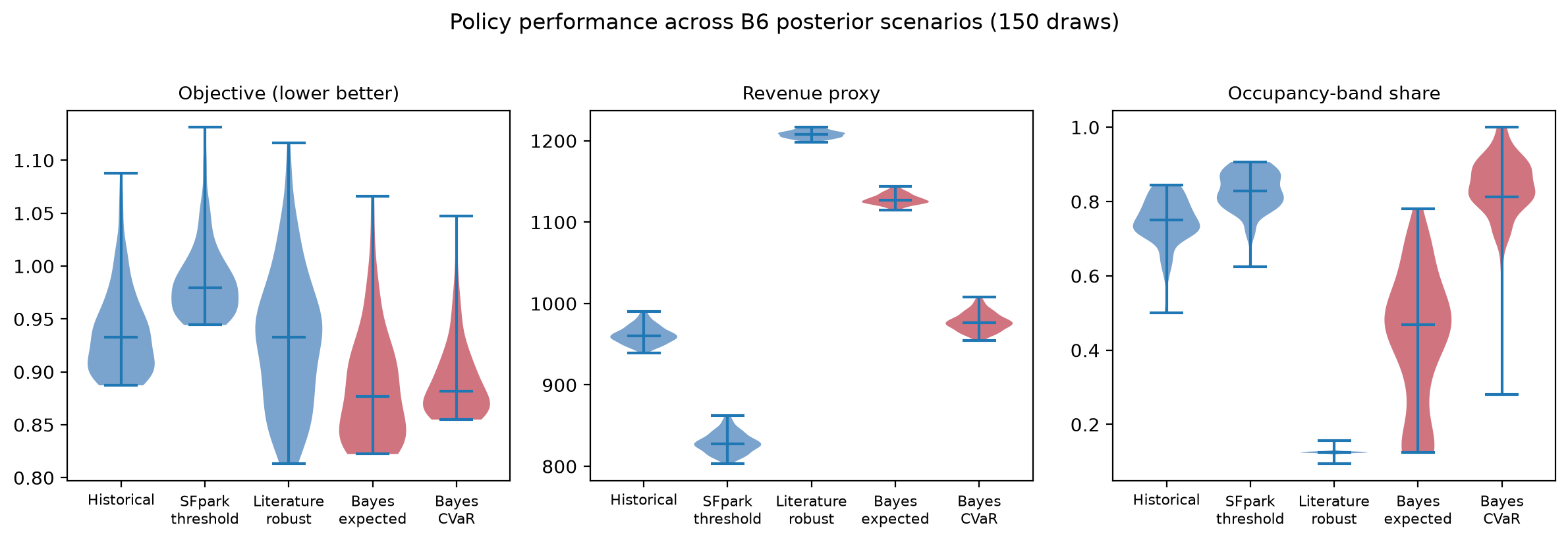}
\caption{Policy performance across 150 posterior scenarios: objective (lower is better),
revenue proxy, and occupancy-band share. The two posterior-based policies are shown in red.}
\label{fig:postpolicy}
\end{figure}

We evaluate five fixed tariffs on 150 held-out joint posterior draws: the historical
tariff, the SFpark threshold rule, the literature-calibrated robust CVaR policy, and the
two posterior-based policies of~\eqref{eq:postobj}. Figure~\ref{fig:postpolicy} shows the
per-draw distributions and Table~\ref{tab:postclaims} the main posterior probabilities.

Among the optimized tariffs, only the posterior CVaR policy keeps occupancy-band share at
the level of the threshold rule (mean $0.83$ for both). At the same time it attains a lower
objective than the threshold rule and the historical tariff in every one of the 150 draws,
and higher revenue than both. Relative to the literature-calibrated robust benchmark, which
pursues revenue under an elasticity ambiguity set considerably wider than the posterior,
the posterior CVaR policy has a lower objective with posterior probability $0.85$ and
higher band share with probability $1.00$, at the cost of lower revenue. The expected-loss
variant attains the lowest mean objective of all policies, but its band-share distribution
is wide; the CVaR variant gives up a small amount of expected objective in exchange for
much more stable band performance. A $95\%$ revenue floor relative to the historical tariff
is satisfied with posterior probability $1.00$.

\begin{table}[!t]
\caption{Posterior policy probabilities on 150 held-out posterior draws
(objective lower is better).}
\label{tab:postclaims}
\centering
\scriptsize
\begin{tabular}{@{}lr@{}}
\toprule
Claim & Posterior probability \\
\midrule
$P(J_{\mathrm{CVaR}}<J_{\mathrm{threshold}})$ & $1.00$ \\
$P(J_{\mathrm{CVaR}}<J_{\mathrm{historical}})$ & $1.00$ \\
$P(J_{\mathrm{CVaR}}<J_{\mathrm{lit.\,robust}})$ & $0.85$ \\
$P(J_{\mathrm{expected}}<J_{\mathrm{lit.\,robust}})$ & $0.97$ \\
$P(R_{\mathrm{CVaR}}>R_{\mathrm{threshold}})$ & $1.00$ \\
$P(R_{\mathrm{CVaR}}>R_{\mathrm{historical}})$ & $1.00$ \\
$P(\mathrm{band}_{\mathrm{CVaR}}>\mathrm{band}_{\mathrm{lit.\,robust}})$ & $1.00$ \\
$P(R_{\mathrm{CVaR}}\ge 0.95\,R_{\mathrm{historical}})$ & $1.00$ \\
\bottomrule
\end{tabular}
\end{table}

\subsection{Multi-district calendar evidence}
Table~\ref{tab:multi} summarizes 20 calendar windows in each of
three districts. Qualitative results persist everywhere: threshold leads band share; coupled
NLP raises revenue by \textbf{21--28\%} and reduces neighbor disparity by about
\textbf{80--83\%} relative to threshold. Window-level paired revenue gains have strictly
positive 95\% bootstrap confidence intervals in all three districts (e.g., Fillmore
$[490,629]$; Mission $[716,942]$; Marina $[361,441]$). Band-share gaps are large but stable
(district-level standard deviations of threshold band share are $0.13$--$0.22$). The
consistency across districts suggests that the Fillmore results are not specific to that
district.

\begin{table}[!t]
\caption{Multi-district calendar backtest (12 blocks, horizon 6, 20 windows each).
$\Delta$rev\% and $\Delta$disp\% are coupled vs threshold.}
\label{tab:multi}
\centering
\scriptsize
\begin{tabular}{@{}lrrrrrr@{}}
\toprule
District & Windows & Thr band & Coupled band & $\Delta$rev\% & $\Delta$disp\% & Rev.\ CI$_{95}$ \\
\midrule
Fillmore & 20 & 0.843 & 0.148 & $+28.0$ & $-79.6$ & $[490,629]$ \\
Mission  & 20 & 0.918 & 0.126 & $+28.1$ & $-83.5$ & $[716,942]$ \\
Marina   & 20 & 0.834 & 0.130 & $+20.9$ & $-82.6$ & $[361,441]$ \\
\bottomrule
\end{tabular}
\end{table}

On Fillmore's 24-window study that includes robust CVaR, threshold band share is $0.846$;
robust CVaR leads revenue ($2564$) and disparity ($0.018$) among optimized policies, with
$120/120$ successful solves.

\subsection{Solution quality}
On $3\times 3$, $4\times 3$, and $4\times 4$ instances, Ipopt's best multistart objective
improves on the best feasible point among $447$--$497$ increment-grid candidates by
$7$--$14\%$. On the $12\times 6$ instance, all seven starts agree to $10^{-8}$ in objective.
Local-optima concerns remain formally open for the full robust NLP, but dispersion evidence
on the deterministic plant is tight.

\section{Discussion and limitations}
\label{sec:discuss}
The main results are stable across the tested elasticity ranges, graph constructions, and districts. If band
compliance is the sole priority under the plant, threshold rules remain strong; the
posterior CVaR policy, however, matches their band share while achieving a lower objective
and higher revenue in every posterior draw. Literature ambiguity does not overturn the
revenue and spatial advantages inside $\varepsilon\in[-0.6,-0.1]$, and replacing that
ambiguity set with the estimated posterior narrows the scenario distribution enough that
band share no longer has to be given up in exchange for tail protection. The
epsilon-constraint frontiers make the same point at the level of the efficient set: under
the posterior, moderate revenue floors cost more band performance than the assumed
calibrations suggest, the knee arrives earlier, and very high revenue floors cost less, so
an agency choosing a floor from the assumed-elasticity frontier would misjudge the tradeoff
in both directions.

Several caveats apply to the Bayesian calibration. The negative-support parameterization is
a structural restriction that we impose and justify on economic grounds; it is not
something the data reveal. The posterior combines the SFpark likelihood with an informative
literature prior and partial pooling, and we report prior-to-posterior updating, a
sign-only companion fit, posterior predictive coverage, and leave-one-wave-out stability so
that the influence of the prior can be assessed. The estimate should be read as a posterior
distribution of short-run price response for use in decision-making, not as a causally
identified elasticity, and we have kept to that language throughout. The aggregation level
is empirically important: the coarser block-wave model failed posterior predictive checks,
whereas the pricing-band model did not. For these data, the aggregation level matters as
much as the model specification.

Limits remain. True geocoded walking networks are unavailable in the ingested files. Revenue
and cruising are proxies. Ipopt is local. Counterfactuals are plant-based. Results are
neighborhood-scale. Deployment would require operational validation, enforcement data, and
stakeholder weight elicitation beyond this study.

\section{Conclusion}
\label{sec:concl}
We estimated a hierarchical posterior distribution of parking-price response from
disaggregated SFpark pricing-band observations and propagated that posterior through a
spatially coupled constrained pricing model. Relative to the earlier design, an externally
assumed elasticity is replaced by an internally estimated posterior ($\varepsilon=-0.241$,
$95\%$ CrI $[-0.325,-0.176]$, with district heterogeneity), and a literature ambiguity set
is replaced by posterior scenarios. Policy performance can then be stated
probabilistically: the posterior CVaR tariff attains a lower operational objective than the
SFpark threshold rule and the historical tariff with posterior probability one, matches
threshold-rule band compliance, and satisfies a $95\%$ revenue floor with posterior
probability one. Demand uncertainty can therefore be incorporated directly into the pricing problem rather than handled through local threshold rules.

\section*{Acknowledgement}
We gratefully acknowledge Dr. Neeta Maitre, department of Computer Engineering, Cummins College for her guidance and continuous support during this research. We also thank Cummins College of Engineering for Women Pune for providing the necessary resources and academic environment that made the completion of this work possible.

SFpark data are products of the San Francisco Municipal Transportation Agency and are used
under SFMTA terms; they are not redistributed with this manuscript. Computations use Ipopt
and Pyomo.

\end{document}